
\input amstex.tex
\documentstyle{amsppt}

\footline={\hss{\vbox to 2cm{\vfil\hbox{\rm\folio}}}\hss}
\nopagenumbers

\def\txt#1{{\textstyle{#1}}}
\baselineskip=13pt
\def\hf{{\textstyle{1\over2}}}
\def\a{\alpha}\def\b{\beta}
\def\d{{\,\roman d}}
\def\e{\varepsilon}

\def\G{\Gamma}
\def\k{\kappa}
\def\s{\sigma}

\def\={\;=\;}

\def\zt{\zeta(\hf+it)}
\def\zr{\zeta(\hf+ir)}

\def\no{\noindent}  
\def\R{\Re{\roman e}\,} \def\I{\Im{\roman m}\,} \def\s{\sigma}
\def\z{\zeta}

\def\no{\noindent} 

\def\H{H_j^3({\txt{1\over2}})}  \def\={\,=\,}
\def\hf{{\textstyle{1\over2}}}
\def\txt#1{{\textstyle{#1}}}

\def\z#1{|\zeta (\hf +i#1)|^4}

\font\tenmsb=msbm10
\font\sevenmsb=msbm7
\font\fivemsb=msbm5
\newfam\msbfam
\textfont\msbfam=\tenmsb
\scriptfont\msbfam=\sevenmsb
\scriptscriptfont\msbfam=\fivemsb
\def\Bbb#1{{\fam\msbfam #1}}

\def \NN {\Bbb N}

\def \RR {\Bbb R}

\font\cc=cmcsc10 at 8pt

\font\ff=cmr8
\def\txt#1{{\textstyle{#1}}}
\baselineskip=13pt
\def\hf{{\textstyle{1\over2}}}
\def\a{\alpha}\def\b{\beta}
\def\d{{\,\roman d}}
\def\e{\varepsilon}

\def\G{\Gamma}
\def\k{\kappa}
\def\s{\sigma}

\def\={\;=\;}

\def\zt{\zeta(\hf+it)}
\def\zr{\zeta(\hf+ir)}

\def\no{\noindent}  
\def\R{\Re{\roman e}\,} \def\I{\Im{\roman m}\,} \def\s{\sigma}
\def\z{\zeta}

\def\d{\,{\roman d}} 
\def\R{\Re{\roman e}\,} \def\I{\Im{\roman m}\,} \def\s{\sigma}

\def\H{H_j^3({\txt{1\over2}})}  \def\={\,=\,}
\font\teneufm=eufm10
\font\seveneufm=eufm7
\font\fiveeufm=eufm5
\newfam\eufmfam
\textfont\eufmfam=\teneufm
\scriptfont\eufmfam=\seveneufm
\scriptscriptfont\eufmfam=\fiveeufm
\def\mathfrak#1{{\fam\eufmfam\relax#1}}

\font\tenmsb=msbm10
\font\sevenmsb=msbm7
\font\fivemsb=msbm5
\newfam\msbfam
     \textfont\msbfam=\tenmsb
      \scriptfont\msbfam=\sevenmsb
      \scriptscriptfont\msbfam=\fivemsb
\def\Bbb#1{{\fam\msbfam #1}}

\def \NN {\Bbb N}

\def \RR {\Bbb R}

  \def\rightheadline{{\hfil{\ff
  On sums of Hecke series in short intervals}\hfil\tenrm\folio}}

  \def\leftheadline{{\tenrm\folio\hfil{\ff
  Aleksandar Ivi\'c }\hfil}}
  \def\emptyheadline{\hfil}
  \headline{\ifnum\pageno=1 \emptyheadline\else
  \ifodd\pageno \rightheadline \else \leftheadline\fi\fi}

  \topmatter
\title
ON SUMS OF HECKE SERIES IN SHORT INTERVALS
\endtitle
\author   Aleksandar Ivi\'c \endauthor
\address{
Aleksandar Ivi\'c, Katedra Matematike RGF-a
Universiteta u Beogradu, \DJ u\v sina 7, 11000 Beograd,
Serbia (Yugoslavia).}
\endaddress
\email {\tt aleks\@ivic.matf.bg.ac.yu,
aivic\@rgf.rgf.bg.ac.yu} \endemail
\abstract
{\ff We have
$$
\sum_{K-G\le\k_j\le K+G} \a_j\H \;\ll_\e\; GK^{1+\e}
$$
for
$K^{\e}  \;\le\; G \;\le \; K$, where
$\a_j = |\rho_j(1)|^2(\cosh\pi\kappa_j)^{-1}$, and
$\rho_j(1)$ is the first Fourier coefficient of the Maass wave form
corresponding to the eigenvalue $\lambda_j = \k_j^2 + {1\over4}$
to which the Hecke series $H_j(s)$ is attached. This result yields
the new bound $H_j(\hf) \;\ll_\e\; \k_j^{{1\over3}+\e}.$}

\medskip\no
{\cc R\'esum\'e}.
{\ff On a
$$
\sum_{K-G\le\k_j\le K+G} \a_j\H \;\ll_\e\; GK^{1+\e}
$$
pour
$K^{\e}  \;\le\; G \;\le \; K$, ou
$\a_j = |\rho_j(1)|^2(\cosh\pi\kappa_j)^{-1}$, et
$\rho_j(1)$ est le premier coefficient de Fourier de forme de Maass
correspondant \`a la valeur propre  $\lambda_j = \k_j^2 + {1\over4}$
 \`a laquelle le s\'erie de Hecke  $H_j(s)$ est attach\'ee.
Ce r\'esultat fournit  l'estimation nouvelle
$H_j(\hf) \;\ll_\e\; \k_j^{{1\over3}+\e}.$}

\endabstract
\endtopmatter
\title par Aleksandar Ivi\'c\endtitle

\def\DJ{\leavevmode\setbox0=\hbox{D}\kern0pt\rlap
{\kern.04em\raise.188\ht0\hbox{-}}D}

\head 1. Introduction and statement of results
\endhead

\bigskip\medskip
The purpose of this paper is to obtain a bound for sums of Hecke
series in short intervals which, as a by-product, gives a new bound
for $H_j(\hf)$. We begin by stating
briefly  the necessary notation and some
results involving the spectral theory of the non-Euclidean Laplacian.
For a competent and extensive account of
spectral theory the reader is referred to Y. Motohashi's monograph [13].

\medskip
Let $\,\{\lambda_j = \kappa_j^2 + {1\over4}\} \,\cup\, \{0\}\,$ be the
eigenvalues (discrete spectrum) of the hyperbolic Laplacian
$$
\Delta=-y^2\left({\left({\partial\over\partial x}\right)}^2 +
{\left({\partial\over\partial y}\right)}^2\right)
$$
acting over the Hilbert space composed of all
$\Gamma$-automorphic functions which are square integrable with
respect to the hyperbolic measure. Let $\{\psi_j\}$ be a maximal
orthonormal system such that
$\Delta\psi_j=\lambda_j\psi_j$ for each $j\ge1$ and
$T(n)\psi_j=t_j(n)\psi_j$  for each integer $n\in\NN$, where
$$
\bigl(T(n)f\bigr)(z)\;=\;{1\over\sqrt{n}}\sum_{ad=n}
\,\sum_{b=1}^df\left({az+b\over d}\right)
$$
is the Hecke operator. We shall further assume that
$\psi_j(-\bar{z})=\e_j\psi_j(z)$ with $\e_j=\pm1$. We
then define ($s = \s + it$ will denote a complex variable)
$$
H_j(s)\;=\;\sum_{n=1}^\infty t_j(n)n^{-s}\qquad(\s > 1),
$$
which we call the Hecke series associated with the Maass wave form
$\psi_j$(z), and which can be continued to an entire function.
The Hecke series satisfies the functional equation
$$
H_j(s) = 2^{2s-1}\pi^{2s-2}\G(1-s+i\k_j)\G(1-s-i\k_j)
(\e_j\cosh(\pi\k_j)-\cos(\pi s))H_j(1-s),
$$
which by the Phragm\'en--Lindel\"of principle (convexity) implies the bound
$$
H_j(\hf) \;\ll_\e \; \k_j^{{1\over2}+\e}.\leqno(1.1)
$$

\medskip\no
It is known that $H_j(\hf) \ge 0$ (see Katok--Sarnak
[8] and for the proofs of (1.2)--(1.4) see [11] or [13]),
and
$$
\sum_{\k_j\le K}\a_jH_j^2(\hf)
\= (A\log K  + B)K^2 + O(K\log^6K)\qquad(A > 0).\leqno(1.2)
$$
Here as usual  we put
$$
\a_j = |\rho_j(1)|^2(\cosh\pi\kappa_j)^{-1},
$$
where $\rho_j(1)$ is the first Fourier coefficient of  $\psi_j(z)$.
Moreover we have
$$
\sum_{\k_j\le K}\a_jH_j^4(\hf)  \ll K^2\log^{15}K\leqno(1.3)
$$
and
$$
\sum_{j=1}^\infty \a_jH_j^3(\hf)h_0(\k_j) =
\left({8\over3} + O\left({1\over\log K}\right)\right)\pi^{-3/2} K^3G
\log^3K\leqno(1.4)
$$
with
$$
K^{{1\over2}}\log^5K \le G \le K^{1-\e}, \leqno(1.5)
$$
$$
h_0(r)
= (r^2 + {\txt{1\over4}})\left(e^{-\left({r-K\over G}\right)^2}
+ e^{-\left({r+K\over G}\right)^2}\right).\leqno(1.6)
$$
Apart from its intrinsic interest, the asymptotic formula (1.4) has an
important application in the theory of the Riemann zeta-function.
Namely it immediately implies
that there are infinitely many $\k$ such that
$$
\sum_{\k_j=\k}\a_j\H \;>\;0,
$$
which is essential in   establishing $\Omega$--results for the function
$E_2(T)$, which represents the error term in the asymptotic formula
for the fourth moment of $|\zt|$ (see [13, Chapter 5]).
Instead of the sum in (1.4) we shall consider the sum
$\sum\limits_{K-G\le \kappa_j \le K+G}\a_j\H$ and
seek an upper bound for it, which is especially
interesting when $G = K^\e$. In that  case it follows
from (1.1) and (1.2) (or from (1.3), or  from (1.4))  that
$$
\sum_{K-K^\e\le \kappa_j
\le K+K^\e}\a_j\H \ll_\e K^{{3\over2}+\e},\leqno(1.7)
$$
where here and later $\e > 0$ denotes arbitrarily small constants,
not necessarily the same ones at each occurrence.  We can suppose that
$$
\sum_{K-K^\e\le \kappa_j
\le K+K^\e}\a_j\H \ll_\e K^{1+\a+\e}\qquad(0 \le \a \le \hf),\leqno(1.8)
$$
and it is reasonable to expect that (1.8) holds with  $\a = 0$.
This is indeed so, and is the content of the following

\medskip
THEOREM.
{\it We have}
$$
\sum_{K-G\le\k_j\le K+G} \a_j\H \;\ll_\e\; GK^{1+\e}\leqno(1.9)
$$
{\it for}
$$
K^{\e}  \;\le\; G \;\le \; K.\leqno(1.10)
$$

\medskip\no
In view of the convention made above on the use of $\e$'s, the
above result strictly speaking means that, for given $\e$ sufficiently small,
the bound (1.9) holds with $GK^{1+\e_1}$ and $\lim\limits_{\e\to0}\e_1=0$,
provided that (1.10) holds.

\medskip
{\bf Corollary 1}. We have (1.8) with $\a = 0\,$.

\medskip
From  $H_j(\hf) \ge 0$ and the bound
$$
\a_j \= {|\rho_j(1)|^2\over\cosh(\pi\k_j)} \;\gg_\e\; \k_j^{-\e}
$$
of H. Iwaniec [6] we obtain

\medskip
{\bf Corollary 2}.
$$
H_j(\hf) \;\ll_\e\; \k_j^{{1\over3}+\e}.
\leqno(1.11)
$$

\smallskip\no
This seems to be the first unconditional improvement
over (1.1), and represents the limit of our method.
Note that H. Iwaniec [7] obtained (1.11) assuming a certain hypothesis
(the referee remarked that, using a trickier amplifier based
on the equality $\lambda_f(p)^2-\lambda_f(p^2) =1$, Iwaniec observed that
his method actually gives unconditionally
$H_j(\hf) \;\ll_\e\; \k_j^{{5\over12}+\e}$, but this  result
sharper than (1.1) does not seem to have appeared in print).
His paper contains several other interesting results, including
a bound for sums of squares of $H_j(s)$ over $\k_j$'s in short intervals.

\medskip
We remark that W. Luo [10] proved the bound
$$
H_j(\hf + i\k_j) \ll_\e \k_j^{{1\over4}+\e}
$$
by exploiting some special properties of the Hecke series at the points
$s = \hf \pm i\k_j$, but our method certainly cannot give such a sharp
bound for $H_j(\hf)$, for which one expects the bound
$H_j(\hf) \ll_\e \k_j^{\e}$, and more generally one conjectures that
$H_j(\hf+it) \ll_\e (|t|\k_j)^{\e}$. This bound may be viewed as a sort
of the  ``Lindel\"of hypothesis"  for $H_j(\hf)$. Since $H_j(s)$ bears
several analogies (i.e., the functional equation) to $\z^2(s)$,
then the bound (1.11) represents the analogue of the
classical estimate $\zt \ll |t|^{1/6}$.

\smallskip
Cubic moments of automorphic $L$-functions $L_f(s,\chi)$ have
been recently investigated by J.B. Conrey and H. Iwaniec [1].
Although they also exploit the idea of the nonnegativity of
cubes of central values of automorphic $L$--functions,
their methods are quite different from
ours. One of their main results is the bound
$$
\sum_{f\in F^\star }L^3_f(\hf,\chi) \ll_\e q^{1+\e},
$$
where $F^\star$ is the set of all primitive cusp forms of weight $k$
(an even integer $\ge 12$) and level dividing $q$, where $\chi(n)
= ({n\over q})$ for odd, squarefree $q$.

\medskip
{\bf Acknowledgement.} I am very grateful to Prof. Matti Jutila for most
valuable remarks.

\bigskip
\head
 2. Beginning of proof
\endhead
\bigskip

Before we begin the proof, some further notation will be necessary.
If one denotes the left-hand side of (1.4) by ${\Cal C}(K,G)$, then with
$\lambda = C\log K\;(C > 0)$ one has ([13, (3.4.18)], with the
extraneous factor $(1 - (\k_j/K)^2)^\nu$ omitted)
$$\eqalign{
{\Cal C}(K,G) &= \sum_{f\le3K}f^{-{1\over2}}\exp\Bigl(-{\bigl({f\over K}
\bigr)}^\lambda\Bigr){\Cal H}(f;h_0)\cr&
- \sum_{\nu=0}^{N_1}\,\sum_{f\le3K}f^{-{1\over2}}U_\nu(fK)
{\Cal H}(f;h_\nu) + O(1),\cr}                    \leqno(2.1)
$$
with ($h_0(r)$ is given by (1.6))
$$
h_\nu(r) \= h_0(r)\left(1 - \left({r\over K}\right)^2\right)^\nu,\leqno(2.2)
$$
$$
{\Cal H}(f;h) \= \sum_{\nu=1}^7{\Cal H}_\nu(f;h),
$$
$$
{\Cal H}_1(f;h) \= -2\pi^{-3}i\left\{(\gamma - \log(2\pi{\sqrt f}))
({\hat h})'(\hf) + {\txt{1\over4}}({\hat h})''(\hf)\right\}d(f)
f^{-{1\over2}},
$$
$$
{\Cal H}_2(f;h) \= \pi^{-3}\sum_{m=1}^\infty m^{-{1\over2}}d(m)d(m+f)
\Psi^+({m\over f};h)\quad\Big(d(n) = \sum_{\delta|n}1\Big),
$$
$$
{\Cal H}_3(f;h) \= \pi^{-3}\sum_{m=1}^\infty (m+f)^{-{1\over2}}
d(m)d(m+f)\Psi^-(1+ {m\over f};h),                      \leqno(2.3)
$$
$$
{\Cal H}_4(f;h) \= \pi^{-3}\sum_{m=1}^{f-1}m^{-{1\over2}}d(m)d(f-m)
\Psi^-({m\over f};h),
$$
$$
{\Cal H}_5(f;h) \= -(2\pi^3)^{-1}f^{-{1\over2}}d(f)\Psi^-(1;h),
$$
$$
{\Cal H}_6(f;h) \= -12\pi^{-2}i\s_{-1}(f)f^{1\over2}h'(-\hf i),
$$
$$
{\Cal H}_7(f;h) \= -\pi^{-1}\int_{-\infty}^\infty{|\z(\hf+ir)|^4\over
|\z(1+2ir)|^2}\s_{2ir}(f)f^{-ir}h(r)\d r\quad (\s_a(f) = \sum_{d|f}d^a),
$$
where
$$
{\hat h}(s) \= \int_{-\infty}^\infty rh(r){\G(s+ir)\over\G(1-s+ir)}\d r,
$$
$$
\Psi^+(x;h) \= \int_{(\b)}\G^2(\hf-s)\tan(\pi s){\hat h}(s)x^s\d s,
$$
and
$$
\Psi^-(x;h) \= \int_{(\b)}\G^2(\hf-s){{\hat h}(s)\over
\cos(\pi s)}x^s\d s,
$$
with $-{3\over2} < \b < {\hf}$, $N_1$ is a sufficiently large integer,
$$
U_\nu(x) = {1\over2\pi i\lambda}\int_{(-\lambda^{-1})}
(4\pi^2K^{-2}x)^wu_\nu(w)
\G({w\over\lambda})\d w \ll \left({x\over K^2}\right)^{-{C\over\log K}}
\log^2K,
$$
where $u_\nu(w)$ is a polynomial in $w$ of degree $\le 2N_1$, whose
coefficients are bounded. A prominent feature of Motohashi's explicit
expression for ${\Cal C}(K,G)$ is that it contains series and integrals
with the classical divisor function $d(n)$ only, with no quantities
from spectral theory. Therefore the problem of obtaining an
upper bound for ${\Cal C}(K,G)$ is a problem of classical analytic
number theory.

\medskip
Now we are ready to begin with the proof of our result.
We shall start from the obvious bound
$$
\sum_{K-G\le \k_j\le K + G}\a_j\H \;\ll\; K^{-2}{\Cal C}(K,G)
\qquad(K^\e \le G \le K),
\leqno(2.4)
$$
so that the proof of the Theorem reduces to showing that
$$
{\Cal C}(K,G) \;\ll_\e\; K^{3+\e}
G\qquad (K^{\e}\le G \le K) .  \leqno(2.5)
$$
The delicate machinery
of (2.1)--(2.3) was developed by Motohashi in order to establish the
asymptotic formula (1.4), where special care must be taken in order to
produce the (weak) error term $O(1/\log K)$. To achieve this, Motohashi
assumed the bound $G \ge K^{1\over2}\log^5K$ in (1.5), which immediately
rendered several contributions in (2.1) negligibly small. However, in (2.5)
we are not aiming at an asymptotic formula for ${\Cal C}(K,G)$, but only
at an upper bound. To obtain this we could start from first principles,
but it seemed expedient to utilize the machinery of (2.1)--(2.3).
First of all, by going through the proof of (1.4), it is seen that it is
the term $\nu = 0$ in (2.1) whose contributions should be considered, because
the bound for the $\nu$-th term will be essentially the same as the
bound for the term $\nu = 0$, only it will be multiplied by $(G/K)^\nu$.
We note that the factors $\exp(-(f/K)^\lambda)$ and $U_\nu(fK)$ in
(2.1) can be conveniently removed by partial summation.
Next we follow the analysis carried out in [13, pp. 120 and 128-129] to
show that the contribution of $\nu = 1,3,5,6,7$ in (2.3) to (2.1) will be
$\ll_\e K^{3+\e}G$. Indeed we have
$$
{\Cal H}_1(f;h_0) \ll d(f)f^{-1/2}K^3G\log^2K, \quad {\Cal H}_3(f;h_0)
\ll e^{-C\log^2K}\quad(C > 0)
$$
by [13, (3.4.20)-(3.4.24)], and in view of [13, (3.3.44)]
$$
{\Cal H}_5(f;h_0) \ll d(f)f^{-1/2},\quad {\Cal H}_6(f;h_0)
\ll \s_{-1}(f)f^{1/2}K.
$$
Finally to deal with ${\Cal H}_7(f;h_0)$ note that we have
 $1/\z(1+ir) \ll \log (|r| + 1)$, $\zr \ll |r|^{1/6+\e}$ (see [4]) and
$$
\sum_{n=1}^\infty \s_{2ir}(n)n^{-ir-s} \= \z(s-ir)\z(s+ir) \qquad
(r \in \RR, \R s > 1).
$$
Consequently by the Perron inversion formula (see e.g., [4, p. 486])
$$
\sum_{f\le3K}\,\s_{2ir}(f)f^{-{1\over2}-ir} \;\ll_\e\; K^{{1\over3}+\e}
\qquad (K \ll |r| \ll K).
$$
Since the relevant range of $r$ in ${\Cal H}_7(f;h_0)$ is
$|r \pm K| \le G\log K$, it follows that the total contribution of
${\Cal H}_7(f;h_0)$ to (2.1) is $\ll_\e GK^{3+\e}$ if $G$ satisfies (1.10).
Thus it transpires that what is non-trivial is the
contribution to (2.1) of
$$
{\Cal H}_2(f;h_0) \= \pi^{-3}\sum_{m=1}^\infty m^{-{1\over2}}d(m)d(m+f)
\Psi^+({m\over f};h_0),\leqno(2.6)
$$
with $m \le 2f$ (the terms with $m > 2f$ are negligible by [13, (3.4.21)])
and
$$
{\Cal H}_4(f;h_0) \= \pi^{-3}\sum_{m=1}^{f-1}m^{-{1\over2}}d(m)d(f-m)
\Psi^-({m\over f};h_0).\leqno(2.7)
$$

\medskip\no
We begin with the contribution of (2.6) for $m \le 2f$, noting that
by [13, (3.4.20)] we have, for $m \le 2f$ and suitable $c > 0$,
$$
\Psi^+({m\over f};h_\nu) \ll K^3G\left({G\over K}\right)^\nu
\exp\left(-cG^2{m\over f}\right)
+ {f\over m}\exp(-{\txt{1\over4}}\log^2K),\leqno(2.8)
$$
which clearly shows that the contribution of the portion of (2.6)
with $m \le 2f$ is negligibly small if (1.5) holds. Our idea is
to evaluate the relevant integrals arising from $\Psi^\pm(m/f;h_0)$
explicitly and then to estimate the ensuing exponential sums, which will
permit us to obtain (2.5) with $G$ lying outside of the range given by (1.5).
From (2.8) it follows that the nontrivial contribution of (2.6) with
$m \le 2f$ will consist of the subsum
$$
\pi^{-3}\sum_{G^2/\log^2K\le f\le3K}f^{-{1\over2}}
\sum_{m\le fG^{-2}\log^2K}m^{-{1\over2}}\ldots\;,
\,\leqno(2.9)
$$
where the sum over $m$ is non-empty for $G \le \sqrt{3K}\log K$.
Henceforth we suppose that
$$
K^\e \le G \le K^{{1\over2}-\e},\leqno(2.10)
$$ which
is actually sufficient for the proof of the Theorem. Namely for the
range $K^{{1\over2}-\e} \le G \le K^{1-\e}$ the bound (1.9) follows
from (1.4)--(1.5), and for $ K^{1-\e} \le G \le K$ from
$\sum_{\k_j\le K}\a_j\H \ll K^2\log^CK$, with an appropriate
change of $\e$ in (1.9).  Now we shall use the formula
after [13, (3.3.39)] with $x = m/f = o(1)$ (as $K\to\infty$), namely
$$\eqalign{&
\Psi^+(x;h) =\cr& 2\pi\int_{-\infty}^\infty rh(r)\tanh(\pi r)
\R\left\{{\G^2(\hf + ir)\over\G(1+2ir)}F(\hf+ir,\hf+ir;1+2ir;-{1\over x})
x^{-ir}\right\}\d r,\cr}
\leqno(2.11)
$$
where $F$ is the hypergeometric function. We shall apply a classical
quadratic transformation formula (see [9, (9.6.12)]) for the hypergeometric
function. This is
$$
F(\a,\b;2\b;z) \= \left({1+\sqrt{1-z}\over2}\,\right)^{-2\a}
F\left(\a,\a-\b+\hf;\b+\hf;\left({1-\sqrt{1-z}\over1+\sqrt{1-z}}\,\right)^2
\right),            \leqno(2.12)
$$
so that (2.11) will give
$$\eqalign{&
\Psi^+(x;h_0)  = 4\pi{\sqrt{x}\over\sqrt{x}+ \sqrt{1+x}}
\int_{-\infty}^\infty rh_0(r)\tanh(\pi r)\R\times\cr&
\left\{{\G^2(\hf + ir)\over\G(1+2ir)}
\Bigl({\sqrt{x}+ \sqrt{1+x}\over2}\right)^{-2ir}
F\left(\hf+ir,\hf;1+ir;\left({\sqrt{x}-\sqrt{1+{x}}\over
\sqrt{x}+\sqrt{1+{x}}}\Bigr)^2
\right)\right\}\d r.\cr}\leqno(2.13)
$$
From the definition (1.6)  it is seen that the integral in
(2.13) will make a negligible contribution unless $|r+K| \le G\log K$
and $|r-K| \le G\log K$. Since the contributions of both ranges of $r$
are treated
analogously (the presence of two exponentials in (1.6) is necessitated
by the fact that Motohashi's approach requires $h_0(r)$ to be an
even function of $r$), we shall treat only the latter, noting that
$\tanh(\pi r) = 1 + O(e^{-K})$ for $|r-K| \le G\log K$. For $|z| < 1$
one has, by the defining property of the hypergeometric function,
$$\eqalign{
F(\a,\b;\gamma;z) &\= \sum_{k=0}^\infty {(\a)_k(\b)_k\over(\gamma)_kk!}z^k\cr&
\= \sum_{k=0}^\infty{\a(\a+1)\ldots(\a+k-1)\b(\b+1)\ldots(\b+k-1)\over
\gamma(\gamma+1)\ldots(\gamma+k-1)k!}z^k.\cr}\leqno(2.14)
$$
We insert (2.14) in (2.13) with
$\a = \hf + ir,\, \b = \hf,\, \gamma = 1 + ir$,
$$
z = \left({\sqrt{x}-\sqrt{1+{x}}\over
\sqrt{x}+\sqrt{1+{x}}}\,\right)^2
= {(\sqrt{x}+\sqrt{1+{x}}\,)}^{-4} = 1 - 4\sqrt{x} + O(x)  < 1 - 5\sqrt{x},
$$
since $m \le fG^{-2}\log^2K$ yields
$x = m/f = o(1)$. In view of the absolute convergence of the series
in (2.14), the resulting relevant expression in (2.13) will be
$$
{4\pi\sqrt{x}\over\sqrt{x}+ \sqrt{1+x}}\sum_{k=0}^\infty{(\hf)_k\over k!}
{(\sqrt{x}+ \sqrt{1+x}\,)}^{-4k}\R I_k,\leqno(2.15)
$$
where
$$
I_k = \int\limits_{K-G\log K}^
{K+G\log K}r(r^2+{\txt{1\over4}})e^{-({r-K\over G})^2}
{(\hf+ir)_k\over(1+ir)_k}
\left(\sqrt{x}+ \sqrt{1+x}\over2\right)^{-2ir}
{\G^2(\hf + ir)\over\G(1+2ir)}\d r\leqno(2.16)
$$
with $K^\e \le G \le K^{{1\over2}-\e}$.
Note that  $(\a)_0 \equiv 1$ and for $k \ge 1$
$$\left|{(\hf+ir)_k\over(1+ir)_k}\right| \le 1,
\quad {(\hf+ir)_k\over(1+ir)_k} = 1 + O\left({1\over r}\right)
$$
uniformly in $k$. The contribution of $k \ge K^{1/2}\log^2K$ will
be clearly negligible, by trivial estimation of the tails of the
series in (2.15).
The contribution of each $I_k$ will be analogous, hence it will suffice
to consider in detail only the case $k = 0$. Note that
$$
{(\hf)_k\over k!} \= {(2k)!\over2^{2k}(k!)^2} \;\ll\; {1\over\sqrt{k}}
$$
if we use the well-known approximation
$$
k! \= \sqrt{2\pi}k^{k+{1\over2}}\exp(-k + {\vartheta\over12k})
\qquad(0 < \vartheta < 1).
$$
Therefore we obtain
$$
\eqalign{
& {4\pi\sqrt{x}\over\sqrt{x}+ \sqrt{1+x}}\sum_{k=0}^\infty{(\hf)_k\over k!}
{(\sqrt{x}+ \sqrt{1+x}\,)}^{-4k} \cr&\ll
\sqrt{x}\sum_{k=0}^\infty(k+1)^{-1/2}(1-5\sqrt{x})^k\cr&
\ll \sqrt{x}\left(\sum_{k\le x^{-1/2}}(k+1)^{-1/2}
+  \sum_{k=0}^\infty x^{1/4}(1-5\sqrt{x})^k\right)\cr&
\ll \sqrt{x}\left(x^{-1/4} + x^{1/4}x^{-1/2}\right)
\ll x^{1/4}.\cr}\leqno(2.17)
$$
Then the expression in (2.9) becomes, up to a negligible error,
$$\eqalign{&\R\Bigl\{
{4\over\pi^2}\sum_{0\le k\le K^{1\over2}\log^2K}{(\hf)_k\over k!}
\sum_{G^2\log^{-2} K\le f\le 3K}f^{-{1\over2}}\times
\cr&
\sum_{m\le fG^{-2}\log^2K}
m^{-{1\over2}}x^{1\over2}(\sqrt{x}+\sqrt{1+x})^{-4k-{1\over2}} I_k\Bigr\},
\cr}\leqno(2.18)
$$
where $x = m/f \ll K^{-\e}$. Note that the expression containing $x$
in (2.18) can be conveniently removed by partial summation. For
each $k$ the double sum over $m$ and $f$  in (2.18) (without
the expression containing $x$) will be
$\ll_\e GK^{3+\e}$ uniformly in $k$ (the key fact is that
the oscillating factor does not depend on $k$), as will be shown in the
next section. Then using (2.17) (with $x$ different from
$x = m/f$, but certainly $x \ll K^{-\e}$) it follows that the total
contribution of (2.9) is $\ll_\e GK^{3+\e}$, as asserted.

Thus it suffices to estimate the contribution
coming from  $I_0$ in (2.18), and to simplify the gamma-factors
in (2.16) we use Stirling's formula in the form ($t \ge t_0 > 0$)
$$
\G(s) = \sqrt{2\pi}\,t^{\s-{1\over2}}\exp\left(-\hf\pi t + it\log t -it
+ \hf{\pi i}(\s - \hf)\right)\cdot\left(1 + O_\s\left(
t^{-1}\right)\right),\leqno(2.19)
$$
with the understanding that the $O$--term in (2.19)
admits an
asymptotic expansion in terms of negative powers of $t$. Therefore
we may replace the gamma-factors in (2.16) by
$Cr^{-1/2}e^{-2ir\log2}(1 + O(1/r))$,
and then make the change of variable $r = K + Gu$ to obtain that the
relevant contribution to $I_0$ will be a multiple of
$$
I' \;:=\; G\int_{-\log K}^{\log K}(K+Gu)^{1\over2}((K+Gu)^2+{\txt{1\over4}})
e^{-u^2}(\sqrt{x} + \sqrt{1+x}\,)^{-2iK-2iGu}\d u.
$$
We expand the first two expressions in $I'$ in power series, taking
sufficiently many terms so that the error term will, by trivial
estimation, make a negligible contribution. The integrals with the
remaining terms are evaluated by using the formula
$$
\int_{-\infty}^\infty u^je^{Au-u^2}\d u \= P_j(A)e^{{1\over4}A^2}
\qquad(j = 0,1,2,\ldots,\;P_0(A) = \sqrt{\pi}\,),\leqno(2.20)
$$
where $P_j(z)$ is a polynomial in $z$ of degree $j$, which may be
explicitly evaluated by successive differentiation of the formula
$$
\int_{-\infty}^\infty e^{Au-u^2}\d u \= \sqrt{\pi}e^{{1\over4}A^2},
$$
considered as a function of $A$. We
note that in each integral over $\,[-\log K,\,\log K]\,$ we
may replace the interval of integration with $(-\infty,\,\infty)$,
making a negligible error. Then we use (2.20) with
$$
A = -2iG\log(\sqrt{x} + \sqrt{1+x}\,) \ll G\sqrt{{m\over f}}\quad
\left(x \= {m\over f} \= o(1)\right),
$$
so that in view of the summation condition in (2.9) we have $A \ll \log^2K$.
The main contribution to $I'$ will come from the
term $j = 0$ in (2.20). This is
$$\eqalign{&
GK^{5/2}{(\sqrt{x} + \sqrt{1+x}\,)}^{-2iK} \int_{-\infty}^\infty
{(\sqrt{x} + \sqrt{1+x}\,)}^{-2iGu}e^{-u^2}\d u\cr&
= \sqrt{\pi}GK^{5/2}{(\sqrt{x} + \sqrt{1+x}\,)}^{-2iK}
\exp\left(-G^2\log^2(\sqrt{x} + \sqrt{1+x}\,)\right),\cr}
$$
and it is precisely the factor $(\sqrt{x} + \sqrt{1+x}\,)^{-2iK}$
which is taken into consideration in our analysis and is crucial
for the proof of the final result.

\bigskip
\head 3. Estimates of exponential sums
\endhead
\bigskip

Now we shall insert the above expression in (2.18) (omitting summation
over $k$ and disregarding the expression containing $x$, as was
just explained), to obtain that the relevant expression which is to
be estimated is a multiple of
$$\eqalign{&
GK^{5/2}\sum_{G^2/\log^2K\le f\le3K}f^{-{1\over2}}\sum_{m\le fG^{-2}\log^2K}
m^{-{1\over2}}d(m)d(m+f)\times\cr&
\times
\left(\sqrt{m\over f} + \sqrt{1 + {m\over f}}\,\right)^{-2iK}
\exp\left(-G^2\log^2\left(\sqrt{{m\over f}}
+ \sqrt{1+{m\over f}}\,\right)\right).\cr}
\leqno(3.1)
$$
Therefore we have reduced the problem to the estimation of the double
exponential sum appearing in (3.1).
The exponential factor in (3.1), which is
$$
\ll \exp\left(-{CG^2m\over f}\right)\qquad(C > 0),
$$
is harmless,
and can be removed by partial summation, being monotonic in $m$ or $f$.
The first idea that might occur in estimating the sum in (3.1) is to
treat it as $\sum\sum d(m)d(m+f)\ldots\,$, namely as the binary additive
divisor problem weighted with an exponential factor. For this problem the
error term is precisely evaluated and estimated by Y. Motohashi [12], and
various averages of the error term by Y. Motohashi and the author [5].
However, summation over the ``shift" parameter $f$ in (3.1)
is too ``long" for
such formulas to be successfully applied. Other possibilities are to
use estimates involving one- and two-dimensional exponent pairs,
coupled with the Voronoi summation formula (see [2]--[4]), to exploit the
particular properties of the function $d(n)$. However all these approaches
yield values of $G$ in a range not as large as the one in (1.10).

\smallskip
To prove the Theorem we shall
proceed in the following, essentially elementary way.
First we change the order of summation in (3.1), keeping in mind
that $m \le fG^{-2}\log^2K$. Then with the help of Taylor's
formula we replace $f^{-1/2}$ by $(m+f)^{-1/2}$, taking sufficiently
many terms so that the contribution made by trivial estimation of the
error term is negligibly small. The contribution of the
term $(m+f)^{-1/2}$ will
be dominant. We replace $m+f$ by $n$, use partial summation to
remove the factor $\exp(-G^2\ldots\,)$, and let $m$ and $n$ lie
in $O(\log^2 K)$ subsums where $M < m \le M_1 \le 2M$, $N < n \le N_1 \le 2N$.
Then we are led to the estimation of the expression
$$\eqalign{&
GK^{5/2}\sum_{M < m \le M_1 \le 2M}\,\sum_{N < n \le N_1 \le 2N}
d(m)m^{-{1\over2}}d(n)n^{-{1\over2}}\times\cr&\hskip14mm\times
\exp\left(2iK
\log\left(\sqrt{m\over n-m} + \sqrt{n\over n-m}\,\right)\right),\cr}
\leqno(3.2)
$$
where we may assume that
$$\eqalign{&
K^\e \le G \le K^{{1\over2}-\e},\; M \ll KG^{-2}\log^2K,\;\cr&
MG^2\log^2K \ll N \ll K,\; N \ge K^{1\over2}, \; MN \ge K.\cr}\leqno(3.3)
$$
The first condition in (3.3) is given by
(2.10), and the next two are implied by (3.1).
Further, for $N \le K^{1\over2}$ (keeping in mind that $M \le N$,
because $m \le fG^{-2}\log^2K$)
or for $MN \le K$  we have,  by trivial
estimation, that the contribution of (3.2) is
$$
\ll_\e GK^{5/2+\e}(MN)^{1/2} \ll_\e GK^{3+\e},
$$
as necessary. Next the range of summation over $n$ in (3.2) is
divided into $O(N/N_0)$ subintervals $\Cal I$ of length at most $N_0$,
where $N_0$ is a parameter that will be suitably chosen a little
later, and which satisfies
$$
1 \;\le N_0 \;\le\; N.\leqno(3.4)
$$
Hence the sum to be estimated is
$$
\sum \;:=\; \sum_{M < m \le M_1 \le 2M}\,\sum_{n \in{\Cal I}}
d(m)m^{-{1\over2}}d(n)n^{-{1\over2}}\exp(iF(m,n)),\leqno(3.5)
$$
where
$$ F(m,n) \;:=\;
2K
\log\left(\sqrt{m\over n-m} + \sqrt{n\over n-m}\,\right).
\leqno(3.6)
$$
By the Cauchy-Schwarz inequality we have
$$
\eqalign{
{\left|\sum\right|}^2 &\le \sum_{M < m \le M_1}{d^2(m)\over m}
\sum_{M < m \le M_1}
\left|\sum_{n\in{\Cal I}}d(n)n^{-1/2}e^{iF(m,n)}\right|^2
\cr&
\ll \log^3M\sum_{M < m \le M_1}\sum_{n_1\in{\Cal I}}\sum_{n_2\in{\Cal I}}
d(n_1)d(n_2)(n_1n_2)^{-1/2}e^{i(F(m,n_1)-F(m,n_2))} \cr&
\ll_\e K^\e\left(N_0MN^{-1} + \sum_{n_1\not=n_2\in{\Cal I}}
N^{-1}\left|\sum_{M < m \le M_1}e^{i(F(m,n_1)-F(m,n_2))}\right|\right).\cr}
$$
If $F_m(m,n)$ denotes the partial derivative of $F(m,n)$
with respect to $m$, then
$$
F_m(m,n) = {Kn\over (n-m)\sqrt{mn}},
$$
and we obtain
$$
\left|F_m(m,n_1) - F_m(m,n_2)\right| \;\asymp\; KM^{-1/2}N^{-3/2}|n_1 - n_2|.
$$
By hypothesis  $|n_1 - n_2| \le N_0$,  thus we  have
$$
\left|F_m(m,n_1) - F_m(m,n_2)\right| \le \hf
$$
if with suitable $C > 0$ we choose
$$
N_0 \;=\; CN^{3/2}M^{1/2}K^{-1}.\leqno(3.7)
$$
Therefore by standard estimates (see e.g.,
[4, Lemma 1.2 and Lemma 2.1]) we have
$$
\sum_{M < m \le M_1}e^{i(F(m,n_1)-F(m,n_2))} \;\ll\;
{M^{1/2}N^{3/2}\over K|n_1-n_2|}\qquad(n_1 \not = n_2, n_1 \in {\Cal I},
n_2 \in {\Cal I}).
\leqno(3.8)
$$
Hence by using (3.8) we obtain
$$
{\left|\sum\right|}^2 \ll_\e K^\e(MN_0N^{-1} + M^{1/2}N_0N^{1/2}K^{-1}).
$$
Consequently the contribution of (3.5) will be, since $M \le N$ by (3.3),
$$\eqalign{&
\ll_\e GK^{{5/2}+\e}NN_0^{-1}(M^{1/2}N_0^{1/2}N^{-1/2} +
N_0^{1/2}M^{1/4}N^{1/4}K^{-1/2})\cr&
\ll_\e GK^{3+\e}(M/N)^{1/4} + GK^{{5/2}+\e}N^{1/2} \ll_\e GK^{3+\e}.\cr}
$$
It remains to check that $N_0$, given by (3.7), verifies (3.4). We have
$$
N_0 = CN^{3/2}M^{1/2}K^{-1} \le N
$$
for $CN^{1/2}M^{1/2} \le K$, which is true in view of $M \le N \le K$.
Also   $NM \ge K$ may be assumed in view of (3.3), and therefore
$$
N_0 = CN^{3/2}M^{1/2}K^{-1} \ge 1
$$
holds for $C(NM)^{1/2}NK^{-1} \gg NK^{-1/2} \ge 1$, that is for,
$N \;\ge\; K^{1/2}$, which is again true by (3.3). Thus
the contribution of (3.1) is $\ll_\e GK^{3+\e}$, and
consequently the total contribution of ${\Cal  H}_2(f;h)$ is also
$\ll_\e GK^{3+\e}$, as asserted.

\bigskip
\head 4. Completion of proof\endhead

\bigskip\no
To finish the proof we have yet to deal with the sum in (2.7). For
$0 < x < 1$ and $-{3\over2} < \b < -\hf$ we have [13, (3.3.45)]
$$\eqalign{&
\Psi^-(x;h) = \cr& = \,\int\limits_0^\infty\Biggl\{\int\limits_
{(\b)}x^s(y(y+1))^{s-1}
{\G^2(\hf-s)\over\G(1-2s)\cos(\pi s)}\d s\Biggr\}
\Biggl\{\int\limits_{-\infty}^\infty rh(r)\left({y\over y+1}\right)^{ir}\d r
\Biggr\}\d y,\cr}\leqno(4.1)
$$
with $h(r)$ given by (1.6) and (2.2). Similarly as in the analysis concerning
(2.11), for $h = h_0$, we may consider only
the ranges $|r+K| \le G\log K$ and
$|r-K| \le G\log K$, and we turn our attention to the latter.
Namely for $|r\pm K| \ge G\log K$  we interchange the order of integration
and in the $y$-integral we integrate the subintegral over $(0,1]$ by parts
to obtain that the contribution is $\ll x^\b\exp(-\hf\log^2K)$.
Therefore the dominant contribution of the $r$-integral will be
$$
\eqalign{&
\int_{K-G\log K}^{K+G\log K}r(r^2 + {\txt{1\over4}})e^{-(r-K)^2G^{-2}}
\left({y\over y +1}\right)^{ir}\d r\cr&
= G\int_{-\log K}^{\log K}(K+Gu)((K+Gu)^2+{\txt{1\over4}}))e^{-u^2}
\left({y\over y +1}\right)^{iK}\left({y\over y +1}\right)^{iuG}\d u.\cr}
$$
We simplify the expression in the first two brackets in the last
integral and use (2.20) with $A = iG\log y/(y+1)$ and $P_1(A) = \hf\sqrt{\pi}
A$. Then the above expression equals $O(\exp(-\hf\log^2K))$ plus
$$\eqalign{&
e^{iK\log{y\over y +1}}e^{-{1\over4}G^2\log^2{y\over y +1}}\left(
\sqrt{\pi}GK^3 + {3\sqrt{\pi}\over2}iG^2K^2\log{y\over y+1}\right)
\cr&+ O\left(KG^3e^{-{1\over8}G^2\log^2{y\over y +1}}\right).\cr}\leqno(4.2)
$$
In view of the exponential factor in (4.2) we may truncate the
$y$--integral in (4.1) at $G/\log K$ with a negligible error. Therefore
the contribution of the $O$--term in (4.2) is, with $\b = \e - 3/2$,
$$
\ll KG^3\int_{G/\log K}^\infty x^\b
y^{2\b-2}\d y \ll  \left({m\over f}\right)^{\e-3/2}KG^{\e-1}.
$$
The total contribution of this expression is  $\ll_\e K^{3+\e}G^{-1}$.
The main terms in (4.2) are treated analogously, and it is the first
one which will make a larger contribution, so it will be treated in detail.
The relevant part of $\Psi^-(x;h)$ will be
$$
\sqrt{\pi}GK^3\int\limits_{(\b)}{x^s\G^2(\hf-s)\over\G(1-2s)\cos(\pi s)}
\Bigl(\,\int\limits_{{G\over\log K}}^
\infty(y(y+1))^{s-1}e^{iK\log{y\over y +1}}
e^{-{1\over4}G^2\log^2{y\over y +1}}\d y\Bigr)\d s.\leqno(4.3)
$$
In view of Stirling's formula  and
$$
|\cos(x + iy)| \= \sqrt{\cos^2x + \sinh^2y}\qquad(x \in \RR,\, y \in \RR),
$$
it follows that the contribution of $|\I s| = |t| > \log^2K$ in (4.3)
will be negligibly small. For $s = \b + it \,(-{3\over2} < \b < -\hf,\,
|t| \le \log^2K)$ we write the integral over $y$ in (4.3) as
$$
I \;:=\;  \int_{G/\log K}^\infty(y^2 + y)^{\b-1}e^{iF(y)}
e^{-{1\over4}G^2\log^2{y\over y +1}}\d y\leqno(4.4)
$$
with
$$
F(y)  \;:=\; t(\log y + \log(y+1)) + K\log y - K\log(y+1)
\qquad(|t| \le \log^2K),
$$
so that
$$
F'(y) \= {t\over y} + {t\over y + 1} + {K\over y(y+1)} \;\gg\;
{K\over y^2}
$$
for $y \ll K\log^{-2}K$. We further write
$$
I =  \int_{G/\log K}^{K^{1-\e}} + \int_{K^{1-\e}}^\infty
= I_1 + I_2 = I_1 + O_\e(K^{2\b-1+\e}),
$$
by estimating $I_2$ trivially. In $I_1$ we write
$e^{iF(y)} = (e^{iF(y)})'/(iF'(y))$ and integrate by parts. Note that the
integrated term at $y = G/\log K$ will be negligibly small in view of the
second exponential factor in (4.4), and at $y = K^{1-\e}$ it will be
$\ll_\e K^{2\b-1+\e}$. We obtain
$$\eqalign{
&I_1 = O_\e(K^{2\b-1+\e}) \;- \cr&
- {1\over i}\int\limits_{G\over\log K}^{K^{1-\e}}
\Bigg\{-{F''(y)\over(F'(y))^2}(y^2+y)^{\b-1} + {1\over F'(y)}
(\b-1)(2y+1)(y^2+y)^{\b-2} \cr&
+ {1\over F'(y)}(y^2+y)^{\b-1}\left(-{G^2\over2y(y+1)}\log{y\over y+1}\right)
\Bigg\}\cr&
\quad\times e^{iF(y)}e^{-{1\over4}G^2\log^2{y\over y +1}}\d y.\cr}
$$
The expression in curly brackets is, since in $I_1$ we have
$F'(y) \gg Ky^{-2}$ and $y \gg G/\log K$,
$$
\ll \left({y\over K} + {G^2\over Ky}\right)y^{2\b-2}
\ll {y\over K}\log^2K\cdot y^{2\b-2}.                    \leqno(4.5)
$$
Thus we obtain the same type of exponential integral again, only in place of
the factor $(y^2+y)^{\b-1}$ we obtain an expression whose order
is given by the right-hand side of (4.5). Since $-3 < 2\b < -1$, this
means that if repeat five times integration by parts we shall obtain an
integral which, when majorized, will have a nonnegative exponent of $y$
in the integrand. Trivial estimation of this integral will yield then
$$
I_1 \;\ll_\e\; K^{2\b-1+\e},
$$
and taking $\b = \e - {3\over2}$ we obtain that ($x = m/f$) for $I$ in
(4.4) we have
$$
K^3GI \;\ll_\e\; K^3Gx^\b K^{2\b-1+\e} \;\ll_\e\; G\left({f\over m}
\right)^{3/2}K^{\e-1}.
$$
By using (2.7) we see that this makes a total contribution
of $\ll_\e GK^{1+\e}$ to (2.1), and thus the proof of the  Theorem
is complete.

\vfill\eject\topskip2cm
\bigskip
\Refs
\bigskip
\item{[1]} J.B. Conrey and H. Iwaniec, The cubic moment of central
values of automorphic $L$-functions, Annals Math., in print.

\item{[2]} S.W. Graham and G. Kolesnik, Van der Corput's Method
of Exponential Sums, {\it LMS Lecture Note Series} {\bf 126}, {\it Cambridge
University Press}, Cambridge, 1991.

\item{[3]} M.N. Huxley, Area, Lattice Points, and Exponential Sums,
{\it London Math. Soc. Monographs} {\bf 13}, {\it  Oxford University Press},
Oxford, 1996.

\item {[4]} A. Ivi\'c,  The Riemann zeta-function, {\it John Wiley and
Sons}, New York, 1985.

\item {[5]}   A. Ivi\'c and Y. Motohashi, On some estimates involving
the binary additive divisor problem, {\it Quart. J. Math. (Oxford)}
(2){\bf46}(1995), 471-483.

\item {[6]} H. Iwaniec, Small eigenvalues of Laplacian
for $\G_0(N)$, {\it Acta Arith.} {\bf 56}(1990), 65-82.

\item {[7]} H. Iwaniec, The spectral growth of automorphic
$L$--functions, J. reine angewandte Math.  {\bf428}(1992),  139-159.

\item {[8]} S. Katok and P. Sarnak, Heegner points, cycles and Maass
forms, {\it Israel J. Math.} {\bf84}(1993), 193-227.

\item {[9]} N.N. Lebedev, Special functions and their applications,
{\it Dover Publications, Inc.}, New York, 1972.

\item {[10]} W. Luo, Spectral mean-values of automorphic L-functions
at special points,  in {\it Analytic Number Theory: Proc. of a
Conference in Honor of H. Halberstam, Vol. 2} (eds. B.C. Berndt et al.),
Birkh\"auser, Boston etc., 1996, 621-632.

\item {[11]} Y. Motohashi, Spectral mean values of Maass wave forms,
{\it J. Number Theory} {\bf 42}(1992), 258-284.

\item {[12]} Y. Motohashi, The binary additive divisor problem,
Ann. Sci. l'\'Ecole Norm. Sup. $4^e$ s\'erie {\bf 27}(1994), 529-572.

\item {[13]} Y. Motohashi,  Spectral theory of the Riemann
zeta-function, {\it Cambridge University Press}, Cambridge, 1997.

\bigskip

Aleksandar Ivi\'c

Katedra Matematike RGF-a

Universitet u Beogradu

\DJ u\v sina 7, 11000 Beograd

Serbia and Montenegro

aivic\@rgf.bg.ac,yu, aivic\@matf.bg.ac.yu

\endRefs

\bye